\documentclass[a4paper,reqno]{amsart}
\usepackage{pgf,tikz,pgfplots}
\usepackage{amsmath}
\usepackage{paralist}
\usepackage{graphics}
\usepackage{epsfig}
\usepackage{amsfonts}
\usepackage{amssymb}
\usepackage{hyperref,cite}
\UseRawInputEncoding
\usepackage{amsthm}
\usepackage{enumerate}
\usepackage[mathscr]{eucal}
\usepackage{eqlist}

\def\teich{{Teichm\"{u}ller}}
\def\mobi{{M{\"{o}}bius}}

\newtheorem{theorem}{Theorem}[section]
\newtheorem{lemma}[theorem]{Lemma}

\def\ifl{\iffalse }

\numberwithin{equation}{section}

\numberwithin{equation}{section}

\theoremstyle{remark}

\begin{document}
\baselineskip=1.3\baselineskip
\title[ On Carleson measure]
{On Carleson Measures of  Beltrami Coefficients Being Compatible with Infinitely Generated Fuchsian Groups Related to Denjoy Domian.}

\author{Huo Shengjin}
\address{Department of Mathematics, Tiangong University, Tianjin 300387, China} \email{huoshengjin@tiangong.edu.cn}
%
%
%
%\author{Michel Zinsmeister}
%\address{Universite D'Orleans, MAPMO, Orleans Cedex 2, France} \email{ zins@unvi-orleans.fr}

%\address{College of Mathematics and statistics,Shenzhen University, Shenzhen, 518060, China} \email{hguo@szu.edu.cn}

\thanks{This work was supported by the Science and Technology Development Fund of Tianjin Commission for Higher Education(Grant No.2017KJ095).}
%the National Natural Science Foundation of China (Grant No.11401432 and Grant No.11571172 ).}
%%National Natural Science Foundation of China (Grant No.11401432 and Grant No.11571172 )
%
%\UseRawInputEncoding
%
\begin{abstract}
  Let $\Omega$ be a Carleson-Denjoy domain and $G$ be its covering group. Let $\mu$ be a Beltrami coefficient on the unit disk which is compatible with the  group $G$ . In this paper we show that if
$\displaystyle\frac{|\mu|^{2}}{1-|z|^{2}}dxdy $ satisfies the Carleson condition on the infinite boundary of the Dirichlet fundamental domain of $G$,
then $\displaystyle\frac{|\mu|^{2}}{1-|z|^{2}}dxdy$  is a Carleson measure on the unit disk. We also show that the above property does not hold for Denjoy domian.

\end{abstract}

\maketitle

\section{1 Introduction }

This is a continuous work of our previous paper\cite{H1} which deal with the Beltrami  coefficients being compatible with the finitely generated Fuchsian groups of the second kind. In this paper, we mainly focus our attention on the  Beltrami  coefficients being compatible with the infinitely generated Fuchsian groups.

 We start by reviewing  some basic definitions about Fuchsian groups. In this paper we call a {\mobi} group  $G$ Fuchsian group  if it acts on the unit
disk $\Delta$ of the complex plane $\mathbb{C}$ properly discontinuously and freely.  The limit set of a Fcuhsian group $G$, denoted by $\Lambda(G)$, is the set of accumulation points of the $G$-orbit of any point $z\in\Delta$. A Fuchsian group is said to be of the first kind if its limit set is the entire circle and of the second kind otherwise.  A Fuchsian  group $G$ is of divergence type if $$\Sigma_{g\in G}(1-|g(0)|)=\infty~~\text{or}~~\sum_{g\in G}\exp(-\rho(0,g(0)))=\infty,$$
where $\rho(0,g(0))$ is the hyperbolic distance between $0$ and $g(0)$. Otherwise, we say that it is of convergence type.
All second kind groups are of convergence type. The readers are suggested to see \cite{Be,Da, KMS}  for more details about Fuchsian groups.

 In this paper, we still use $\mathcal{F}_{z}(G)$ to denote the Dirichlet fundamental domain of $G$ centered at $z$. For simplicity, we use the notation $\mathcal{F}$ to express  the Dirichlet fundamental domain $\mathcal{F}_{z}(G)$ of $G$ centered at $z=0.$  We call the intersection of $\overline{\mathcal{F}}$ with the unit circle $\partial\Delta$ the boundary at infinity of $\mathcal{F}$, denoted by $\mathcal{F}(\infty)$.

Recall that a positive measure $\lambda$ defined in a simply connected domain $\Omega$ is called a Carleson measure if there exists some constant C which is independent  of $r$ such that, for all $0<r<diameter(\partial\Omega)$ and $z\in \partial\Omega$,
 $$\lambda(\Omega\cap D(z,r))\leq Cr. $$
 The infimum of all such $C$ is called the Carleson norm of $\lambda,$ denoted by
 $\parallel\lambda\parallel_{*}.$ For more detail about Carleson measure, see \cite{Ga,WZ}.

 Let $G$ be a Fuchsian group and $\mu(z)$ a bounded measurable function on $\Delta$,  then we say $\mu$ is a  $G$-compatible Beltrami coefficient (or complex dilatation) if it satisfies $$||\mu(z)||_{\infty}<1 ~~\text {and} ~~\mu(z)=\mu(g(z))\overline{g'(z)}/g'(z),\eqno (1.1)$$ for every $g\in G$. We use $M(G)$ to denote the set of all $G$-compatible Beltrami coefficients. For a $G$-compatible Beltrami coefficient $\mu$, if the measure
$$\displaystyle\frac{|\mu|^{2}}{1-|z|^{2}}dxdy$$
is a Carleson measure on $\Delta$ and when the Carleson norm is small, then $f_{\mu}(\partial\Delta)$ is a rectifiable (chord-arc) curve, where $f_{\mu}$ is the quasiconformal mapping of the complex plane $\mathbb{C}$ with $i$, $1$ and $-i$ fixed, whose Beltrami coefficient equals to $\mu$ a.e. on the unit disk and equals to zero on the outside of the unit disk. This is essential for the proof of the convergence-type first-kind Fuchsian groups failing to have Bowen's property, see \cite{AZ1}. It is also the method to prove that some convergence-type Fuchsian groups  fail to have Ruelle's property, see \cite{HW, HZ}.

 We say that a measurable function $\mu(z)$ belongs to $ CM^{*}(\Delta)$ if the measure $$\displaystyle\frac{|\mu|^{2}}{1-|z|^{2}}dxdy\in CM(\Delta).\eqno (1.2)$$

The importance of the class $CM^{*}(\Delta)$ lies on the fact that it has wide applications in BMOA-{\teich} space and Weil-Petersson {\teich} space, see \cite{AZ,Bi2, Cu, SW} etc. Hence it is important to investigate under which condition the $G$-compatible Beltrami coefficients belong to $CM^{*}(\Delta).$

It is natural to ask whether or not can we determine a $G$ compatible Beltrami coefficient $\mu$  belonging to  the class $CM^{*}(\Delta)$ by its value on the Dirichlet domains? Quite recently, we have proved

\begin{theorem}\label{main}(\cite{H1})
Let $G$ be a convex cocompact Fuchsian group of the second kind  and $\mathcal{F}$ the Dirichlet  fundamental domain of $G$ centered at $0$. Let $\mu\in M(G)$: if there exists a constant $C$ such that, for any $\xi\in \mathcal{F}(\infty)$(i.e. $\xi$ is in the free edges of $\mathcal{F}$) and for any $0<r<1$,
$$\iint_{B(\xi,r)}\displaystyle\frac{|\mu|^{2}\chi_{\mathcal{F}}}{1-|z|^{2}}dxdy \leq Cr,\eqno (1.3)$$
then $\mu$ is in $CM^{*}(\Delta),$ where $\chi_{\mathcal{F}}$ is the characteristic function of the Dirichlet fundamental domain $\mathcal{F}.$
\end{theorem}

Furthermore, Theorem \ref{main} can be generalized to the finitely generated Fuchsian group of the second kind with some parabolic elements.
We have
\begin{theorem}\label{main1}(\cite{H1})
Let $G$ be a finitely generated Fuchsian group of the second kind with some parabolic elements and $\mathcal{F}$ the Dirichlet fundamental domain of $G$ centered at $0$. Let $\mu\in M(G)$: if there exists a constant $C$ such that, for any $\xi\in \mathcal{F}(\infty)$ and for any $0<r<1$,
$$\iint_{B(\xi,r)}\displaystyle\frac{|\mu|^{2}\chi_{\mathcal{F}}}{1-|z|^{2}}dxdy \leq Cr,\eqno(1.4)$$
then $\mu$ is in $CM^{*}(\Delta).$
\end{theorem}

This theorem means that the Carleson property of the measures which are compatible with the finitely generated Fuchsian groups can be checked from the points  in the set $\mathcal{F}(\infty)$ i.e., the boundary at infinity of the Dirichlet domain $\mathcal{F}$.

Notice that Theorem \ref{main} fails for the case of the finitely generated Fuchsian groups of the first kind (i.e. cocompact groups), since  Bowen \cite{Bo} showed that cocompact groups hold a rigidity property, now called Bowen's property, i.e. the image of the unit circle under any quasiconformal map whose Beltrami coefficient compatible with a cocompact group,  is either a circle or  has Hausdorff dimension bigger than 1. Hence for any $\mu$ being compatible with geometry finite groups, the measure $\displaystyle \frac{|\mu|^{2}}{1-|z|^{2}}dxdy$ is not a Carleson measure.

Recall that a Denjoy domain is a connected open subset $\Omega$ of the extended complex plane $\bar{C}$ such that the complement $E=\bar{C}\setminus\Omega$ is a subset of the real axis $\mathbb{R}.$ In addition, $\Omega$ is called a Carleson-Denjoy domain if there exists a positive constant $C$ such that
 $$|E\bigcap(x-t,x+t)|\geq Ct \eqno (1.5)$$
for all $x\in \mathbb{R}$ and $0,t<\text{diam}(E)$, where $|\cdot|$ denotes the Lebesque measure on $\mathbb{R}$. Let $G$ be a covering group of the unit disk $\Delta$ over $\Omega,$ i.e., a Fuchsian group of  $\Delta$   with quotient $\Delta/G$ conformally equivalent to $\Omega.$  Any two covering groups of $\Omega$ are conjugate, and coversely. It is easy to see that when the boundary of $\Omega$ is  totally disconnected, the covering group $G$ is of infinitely generated Fuchsian group of the first kind and of  convergence type. For such groups $G$, we have

\begin{theorem}\label{main1}
Let $\Omega$ be any Carleson-Denjoy domain with a totally disconnected boundary. Let $G$ be  the covering group of the unit disk $\Delta$ over $\Omega$ and $\mathcal{F}$ the Dirichlet  fundamental domain of $G$ centered at $0$. Let $\mu\in M(G)$: if there exists a constant $C$ such that, for any $\xi\in \mathcal{F}(\infty)$(i.e. $\xi$ is in the free edges of $\mathcal{F}$) and for any $0<r<1$,
$$\iint_{B(\xi,r)}\displaystyle\frac{|\mu|^{2}\chi_{\mathcal{F}}}{1-|z|^{2}}dxdy \leq Cr,\eqno(1.6)$$
then $\mu$ is in $CM^{*}(\Delta),$ where $\chi_{\mathcal{F}}$ is the characteristic function of the Dirichlet fundamental domain $\mathcal{F}.$
\end{theorem}

I thank professor Michel Zinsmeister for pointing out this result to me.

We will show that if a Denjoy domain $\Omega$ does not satisfy the homogeneous property $(1.5)$, i.e. $\Omega$ is not a Carleson-Denjoy domain, Theorem \ref{main1} will not hold. We shall use the rich knowledge about the Ruelle's property about Fuchsian groups to prove the following  result.

\begin{theorem}\label{main2}
  Suppose$(s_n)$  be a sequence of real numbers increasing to infinity and $G$ be the covering group of the surface $S=\mathbb{C}\backslash\{s_n,\,n\geq 0\}.$  There exists a sequence $(s_n)$ and  a constant $C$ such that, for any $\xi \in  \mathcal{F}(\infty)$(i.e. $\xi$ is in the infinity boundary of $\mathcal{F}$) and for any $0<r<1$,
$$\iint_{B(\xi,r)}\displaystyle\frac{|\mu|^{2}\chi_{\mathcal{F}}}{1-|z|^{2}}dxdy \leq Cr,\eqno(1.7)$$
however $\mu$ is not in $CM^{*}(\Delta),$ where  $\mathcal{F}$ is  the Dirichlet fundamental domain of the group $G.$
\end{theorem}

\medskip

 The structure of the rest of the paper is as follows: Section 2, we recall some basic definitions and results on Bowen's property and Ruelle's property which will be used in the proof of Theorem 1.4.

\section{Bowen's property and Ruelle's property revisited}

 Let $G$ be a Fuchsian group and $\mu \in M(G)$. By the measurable Riemann mapping theorem, there exists a unique quasiconformal self-mapping $f^{\mu}$ of $\Delta$ fixing 1,-1 and i, and satisfying
$$\bar{\partial}f^{\mu}=\mu\partial f^{\mu}~~~~~a.e. ~z\in \Delta.$$

Similarly, there exists a unique quasiconformal homeomorphism
$f_{\mu}$ of the plane $\mathbb{C}$ which is holomorphic outside
the unit disk $\Delta$, fixing 1,-1and i, and satisfying

$$\bar{\partial}f_{\mu}=\mu\partial f_{\mu} ~~~~~a.e.~~ z\in \mathbb{D}.$$
Then $f_{\mu}(\partial\Delta)$ is a quasi-circle,
i.e. the image of $\partial \mathbb{D}$ under a quasiconformal
mapping of the plane.

For any $\mu \in M(G)$, $f^{\mu}\circ g\circ f^{-\mu}$ is a {\mobi} transformation for every $g\in G.$ In this case, we call $G'=f^{\mu}\circ G\circ f^{-\mu}$ is a quasiconformal deformation of $G$. The group $G^{*}=f_{\mu}\circ G\circ f_{-\mu}$ is called a quasi-Fuchsian group. We say that a Fuchsian group $G$ has Bowen's property if the limit set of any quasiconformal deformation of $G$ is either a circle or has Hausdorff dimension  $>1$. In 1979, Bowen \cite{Bo} proved that if $G$ is a finitely generated Fuchsian group of first kind without parabolic elements, then the limit set of any quasiconformal deformation of $G$ is either a circle or has Hausdorff dimension  $>1$. Soon, Sullivan \cite{Su1,Su2} extended Bowen's property to all cofinite groups. In 1990, K. Astala and M. Zinsmeister \cite{AZ2} showed that Bowen's property fails for all convergence groups of the first kind. At last, in 2001, C.J. Bishop \cite{Bi1} showed a excellent result for all divergence groups as follows:

\begin{lemma}\label{le}\cite{Bi1}
Suppose $G$ is a divergence type Fuchsian group and $G'=f^{\mu}\circ G\circ f^{-\mu}$ is a quasiconformal deformation of $G$. Then either $f_{\mu}(\partial\Delta)$ is a circle or has Hausdorff dimension  $>1$.
\end{lemma}

 Let us recall that a Fuchsian group $G$ has Ruelle's property if, for any  family of Beltrami coefficients $(\mu_{t})\in M(G)$ which is analytic in $t\in\Delta$, the map $t\mapsto HD(\Lambda(G_{\mu_{t}}))$ is  real-analytic in $\Delta$.
In 1982, Ruelle \cite{Ru} showed that all cocompact groups have this property. In 1997, J.W. Anderson and A.C. Rocha \cite{AR} extended this result to finitely-generated Fuchsian groups without parabolic elements.  In \cite{AZ, AZ1}, Astala and Zinsmeister showed that for Fuchsian groups corresponding to Denjoy-Carleson domains or infinite $d$-dimensional "jungle gym" with $d\geq 3$ , Ruelle's property fails.  It is easy to see that the Fuchsian groups
studied in \cite{AZ1}and \cite{Bi1} are of the first kind.
Very recently, the author and Zinsmeister showed that
\begin{lemma}{\label{le2}}\cite{HZ}
All convergence type Fuchsian groups of the first kind fail to have Ruelle's property.
  \end{lemma}

  In \cite{HZ}, we constructed an infinitely generated Fuchsian group which has Ruelle's property. As the author's known, this is the first concrete example of infinitely generated Fuchsian group which has Ruell's property.

\begin{lemma}\label{le3}\cite{HZ}
 There exists a sequence $(s_n)$ of real numbers increasing to infinity such that the Fuchsian group uniformizing $S=\mathbb{C}\backslash\{s_n,\,n\geq 0\}$ has Ruelle's property.
\end{lemma}

Now we give the proofs of Theorem 1.3 and 1.4.
\section{Proof of Theorem \ref{main1}}
 For all $\xi\in \partial\Delta$ and all $0<r<2$, we need to show that there is a constant $C$ which is independent on $\xi$ and $r$ such that
 $$\iint_{B(\xi,r)\cap \mathbb{D}}\frac{{{{\left| {{\mu(w)}} \right|}^2}}}{{1 - \left| w \right|^{2}}}dudv\leq Cr,$$
  where $B(\xi, r)$ is the disk with center $\xi$ and radius $r$.

  Note that
\begin{eqnarray}
\iint_{B(\xi,r)\cap \mathbb{D}}\frac{{{{\left| {{\mu(w)}} \right|}^2}}}{{1 - \left| w \right|^{2}}}dudv
&=& \iint_{g^{-1}(B(\xi,r)\cap \mathbb{D})}\frac{{{{\sum_{g\in G}\left| {{\mu(g(z))}} \right|}^2\chi_{g(\mathcal{F})}}}}{{1 - \left| g(z) \right|^{2}}} dxdy\\
&=& \sum_{g\in G}\iint_{g^{-1}(B(\xi,r))\cap \mathcal{F}}\frac{{{{\left| {{\mu(g(z))}}\right|}^2}}| g'(z)|^{2}}{{1 - \left| g(z) \right|^{2}}} dxdy\\
&=&\sum_{g\in G}\iint_{g^{-1}(B(\xi,r))\cap \mathcal{F}}\frac{\mid\mu(z)\mid^{2}}{1-\left|z\right|^{2}}\left|g'(z)\right|dxdy
\end{eqnarray}

The equality $(3.3)$ holds since $\mu $ is compatible with the group $G.$
By the statement of the theorem, suppose $C$ is the constant such that for any $\zeta \in \mathcal{F}(\infty),$
\begin{eqnarray}\iint_{B(\zeta,r)}\displaystyle\frac{|\mu|^{2}\chi_{\mathcal{F}}}{1-|z|^{2}}dxdy \leq Cr.
\end{eqnarray}
By some calculation or see (Lemma 2.1, \cite{H1}), we know that for any $\xi\in \mathcal{F}$ and any $0<r<2,$ the inequality (3.4) still holds.

Hence the measure $$\iint_{B(\zeta,r)}\displaystyle\frac{|\mu|^{2}\chi_{\mathcal{F}}}{1-|z|^{2}}dxdy$$
is a Carleson on the domain $g^{-1}(B(\xi,r))\cap \mathcal{F}.$
Combine  with (3.3) we have

\begin{eqnarray}
\iint_{B(\xi,r)\cap \mathbb{D}}\frac{{{{\left| {{\mu(w)}} \right|}^2}}}{{1 - \left| w \right|^{2}}}dudv&\leq&C_1\sum_{g\in G}\int_{\partial(g^{-1}(B(\xi,r))\cap \mathcal{F})}\left|g'(z)\right|ds\\
&=&C_1\sum_{g\in G}\int_{\partial(B(\xi,r)\cap g(\mathcal{F}))}ds
\end{eqnarray}

In \cite{FH} Fernandez and Hamilton showed that for Carleson-Denjoy domain $\Omega^{*}$, $\Gamma$ its covering group and $\mathcal{F}_{0}$ the Dirichlet domain of $\Gamma$ with center $0,$
\begin{eqnarray}
\sum_{\gamma\in \Gamma}\text{length}(\partial(\gamma(\mathcal{F})))<\infty,
\end{eqnarray}
 or see \cite{F}. In fact , Carleson \cite{C} for Carleson-Denjoy domain (in \cite{C} it is called homogeneous), the harmonic measure is  absolutely continuous.

  Note that $$\text{length}~~\partial (B(\xi,r)\cap g(\mathcal{F}))<2~\text{length}~(B(\xi,r)\cap\partial g(\mathcal{F})),$$
  combine with (3.6) we prove the theorem.

\section{Proof of Theorem \ref{main2}}
Since hyperbolic area is unchanged under conformal mapping,
 for convenience, we will  first use the upper half plane $\mathbb{H}=\{z: Im(z)>0\}$ as the covering group of the surface $S.$
Let $\mathcal{D}^{*}_{1}$ be the closed disk with diameter $[0, 2]$ and $\mathcal{D}^{*}_{n}$, $n\geq 2$,  the closed disk with diameter $[2^{n-1}, 2^{n}].$
We consider the domain $$\Omega=\mathbb{H}\setminus((\cup_{n\geq1}\mathcal{D}^{*}_{n})
\cup(\cup_{n\geq 1}(-\mathcal{D}^{*}_{n}))).$$

Let $\phi$ be the conformal mapping from $\Omega$ onto $\mathbb{H}$ fixing $0$, $2$ and $\infty.$ We put $z_{0}=0$, and $z_{n}=\phi(2^{n}), \, n\geq1$ and  $z_{n}=\phi(-2^{-n}), \, n\leq-1.$ Let $\sigma_{n}$ be the reflection with respect to $\partial\mathcal{D}^{*}_{n}$ and $\tau(z)=-\bar{z}$.
By Rubel and Ryff's construction \cite{RR} of the covering group of Riemann surface
$S=\mathbb{C}\setminus \{z_{n}\},$ the Fuchsian group $\Gamma$ generated by $\{\tau\circ\sigma_{n}\}^{\infty}_{n=1}$ uniformities the surface $S$, in the sense that $S\simeq \mathbb{H}/\Gamma.$ By the construction of $\Omega$ we can see that $\Gamma$ is of infinitely generated  and of first kind. It is easy to see that $\Gamma$ contains infinitely many parabolic elements. Hence the Dirichlet domain $\mathcal{F}$ of $\Gamma$ contains countably many cusps and the set $\mathcal{F}(\infty)$, i.e., the infinity boundary  of $\Gamma$, contains countably many points, denoted  $\mathcal{F}(\infty)$ by $\{\zeta_{n}\}^{\infty}_{-\infty}.$

Let $B^{*}_{n}=B(\zeta_{n}, 1)\cap \mathcal{F}$  and $B^{*}=\cup B^{*}_{n}$. In the following we will show that the hyperbolic area of $\cup B^{*}_{n}$ is finite. We first give the hyperbolic area of $ B^{*}_{0}.$ Without loss of generality, we may suppose  $\zeta_{0}=0$ and let $C_{-1}$ and $C_{1}$ be the two infinity sides of the Dirichlet domain of the group $\Gamma$ with $0$ as a vertex, respectively.

Suppose $C_{-1}$ be the circle $$(x-r_{-1})^{2}+y^{2}=r_{-1}^{2}$$ and $C_{1}$  the circle
$$(x+r_{1})^{2}+y^{2}=r_{1}^{2}.$$

Then we have
\begin{align}\text{Area} (B^{*}_{0})&=\iint_{B^{*}_{0}} \displaystyle \frac{1}{4v^{2}}dudv\\
&=\int^{1}_{0}dr\int^{\arccos (\displaystyle\frac{-r}{2r_{-1}})}_{\arccos (\displaystyle\frac{r}{2r_{1}})}\displaystyle\frac{1}{4r\sin^{2}\theta}d\theta.\\
&=\int^{1}_{0}(-cot\theta\big|^{-\displaystyle\frac{r}{4r_{-1}}}_{\displaystyle\frac{r}{4r_{1}}})dr\\
&\leq C(\frac{1}{r_{-1}}+\frac{1}{r_{1}}),
\end{align}
where $C$ is a universal constant.
For any  integer $n\in \mathbb{Z}$,  by some calculation as above, we can get the hyperbolic area of $ B^{*}_{n}.$
We have, for any $n\in \mathbb{Z}$,

\begin{align}
\text{Area}(B^{*}_{n})=&\iint_{B^{*}_{n}} \displaystyle \frac{1}{4v^{2}}dudv\\
&=\int^{1}_{0}dr\int^{\arccos (\displaystyle\frac{-r}{2r_{n-1}})}_{\arccos (\displaystyle\frac{r}{2r_{n}})}\displaystyle\frac{1}{4r\sin^{2}\theta}d\theta.\\
&\leq C(\displaystyle\frac{1}{r_{n-1}}+
\displaystyle\frac{1}{r_{n}}),
\end{align}
where the equality (4.6) is  from  the hyperbolic area of the domain in the upper half plane $\mathbb{H}$ being unchanged under the translation $f(z)=az+b$ along the horizontal direction,  $a, b\in \mathbb{R}.$

Note that the conformal mapping  $\phi$ from $\Omega$ onto $\mathbb{H}$ fixing $\infty,$ we have  $z_{n}=\phi(2^{n})$ is comparable to $2^{n}$ for $n>0$ and $z_{n}=\phi(-2^{-n})$ is comparable to $-2^{-n}$ for $n<0$. Furthermore, the radius $r_{n}$ is comparable to $z_{n}$. Combine with (4.5) and by some easy calculation, we know the hyperbolic area of $B^{*}$ is finite.

In the following, we will show that for any $\zeta_{n}$ in $\mathcal{F}(\infty)=\{\zeta_{m}\}^{\infty}_{-\infty}$  and $0<r<1,$

Note that the limit
$$\lim_{r\rightarrow 0}\int^{\arccos (\displaystyle\frac{-r}{2r_{n}})}_{\arccos (\displaystyle\frac{r}{2r_{n-1}})}\displaystyle\frac{1}{4r\sin^{2}\theta}d\theta
=\displaystyle\frac{1}{8}(\displaystyle\frac{1}{r_{n-1}}+
\displaystyle\frac{1}{r_{n}})$$
 and the sequence $(r_{n})$ is increasing to infinity as $n$ tending to infinity.

Suppose $\mu$ is a measurable function on $\mathbb{H}$ with $\parallel \mu\parallel_{\infty}<1.$
by some calculation as (4.5) to (4.7), for any $0<r<1$, for any $\xi \in  \mathcal{F}(\infty)$(i.e. $\xi$ is in the infinity boundary of $\mathcal{F}$) and for any $0<r<1$,

\begin{align}
&\iint_{B(\xi,r)\cap\mathcal{F}}\displaystyle\frac{|\mu|^{2}}{(\text{Im} (z))^{2}}dxdy\\
&\leq\int^{r}_{0}dr\int^{\arccos (\displaystyle\frac{-r}{2r_{n-1}})}_{\arccos (\displaystyle\frac{r}{2r_{n}})}\displaystyle\frac{1}{4r\sin^{2}\theta}d\theta \leq Cr.
\end{align}
Note that for any $0<r<1,$
\begin{align}
\iint_{B(\xi,r)\cap\mathcal{F}}\displaystyle\frac{|\mu|^{2}}{(\text{Im} (z))}dxdy\leq\iint_{B(\xi,r)\cap\mathcal{F}}\displaystyle\frac{|\mu|^{2}}{(\text{Im} (z))^{2}}dxdy.
\end{align}

Define $$\mu_{\mathcal{F}}(z) = \left\{ \begin{gathered}
\mu(z),\hfill z\in B^{*}\cap \mathcal{F};\hfill\\
0,\hfill~~z \in \mathcal{F}\backslash B^{*} ,\hfill\\
\end{gathered}  \right.\eqno(4.11)
$$
and
$$\mu^{*}(z)=\sum_{\gamma\in \Gamma}\mu_{\mathcal{F}}\circ \gamma^{-1}(z)\chi_{\mathcal{F}}\circ \gamma^{-1}(z).\eqno(4.12)$$ It is easy to see that $\mu^{*}(z)\in M(\Gamma).$

Now it's time to replace the upper half plane $\mathbb{H}$ by the unit disk $\Delta$ as the covering surface. Consider the Cayley transformation $\kappa(z)=\displaystyle \frac{z-i}{z+i}$ from the upper half plane $\mathbb{H}$ onto the unit disk $\Delta.$

By the transformation $\kappa$, we can conjugate the Fuchsian group $\Gamma$ to a Fuchsian group $G$ of the unit disk as $G=\kappa\circ\Gamma\circ\kappa^{-1}$ and draw $\mu^{*}$ defined in (4.12) to the unit disk as $\mu_{0}(z)=\mu^{*}\circ\kappa^{-1}(z)(\frac{\overline{(\kappa^{-1})'}}{(\kappa^{-1})'}).$
Combine (4.8), (4.9), (4.10)and the construction of $\mu_{\mathcal{F}}$, we know that there exists a constant $C$ such that, for any $\xi\in \mathcal{F}(\infty)$(i.e. $\xi$ is in the free edges of $\mathcal{F}_{0}$) and for any $0<r<1$,
$$\iint_{B(\xi,r)}\displaystyle\frac{|\mu_{0}|^{2}\chi_{\mathcal{F}_{0}}}{1-|z|^{2}}dxdy \leq Cr,\eqno(4.13)$$
 where $\mathcal{F}_{0}$ denotes the Dirichlet domain of the group $G$ with center $0$.

In the following we show that $\mu_{0}\in CM^{*}(\Delta).$

It is easy to see that $\mu_{0}\in M(G).$ By the construction of the Denjoy and Lemma \ref{le3}  we know the group $G$ has Ruelle's property. By Lemma \ref{le2} we have $\mu$ is divergence type.   Suppose $\mu_{0}\in CM^{*}(\Delta).$ It is known that in this case $\log (f_\mu')$ belongs to the space $BMOA(\Delta)$ with a norm controlled by the above Carleson measure norm. In particular, when the Carleson norm is small then $\partial{f_\mu(\Delta)}$ is a rectifiable (chord-arc) curve (\cite{Po} or \cite{Se}). Hence this is contradict with Lemma \ref{le} which said that for all divergence type Fuchsian groups, if $G'=f^{\mu}\circ G\circ f^{-\mu}$ is a quasiconformal deformation of $G$, the Hausdorff dimension of $f_{\mu}(\partial\Delta)$ is bigger than one.

\section{ Acknowledgements}

It is my pleasure to thank professor Michel Zinsmeister for inviting me to
the University of Orleans as a visiting scholar for one year and for some discussions on topics related to this paper. The author would also like to thank China Scholar Council for life-expenses in Orl\'{e}ans.


\begin{thebibliography}{99}
 \footnotesize

 \bibitem {AR}J.W. Anderson and A.C. Rocha.  {\em Analyticity of Hausdorff dimension of limit sets of Kleinian groups.} Ann. Acad. Sci. Fenn. Math., Vol 22(1997), 349-364.

\bibitem {AZ}K. Astala and M. Zinsmeister.  {\em {\teich} spaces and BMOA.} Math. Ann., Vol 289(1991), 613-625.

\bibitem {AZ1} K. Astala and M. Zinsmeister. {\em Holomorphic families of quasi-Fuchsian groups.} Ergod.Th and Dynam. sys. Vol 14(1994) 207--212.

\bibitem{AZ2} K. Astala and M. Zinsmeister. {\em Rectifiability in Teichm¨¹ller theory. in Topics in Complex Analysis,} Banach Center Publications Vol 31, (1995), 45-52.


\bibitem {Be} A. F. Beardon. {\em The geometry of discrete group.} Springer-Verlag, 1983.

\bibitem{Bi1} C.J. Bishop. {\em Divergence groups have the Bowen property.} Ann. Math., Vol 154(2001), 205-217.

\bibitem{Bi2} C.J. Bishop. {\em Compact deformations of Fuchsian group.} J. D'analyse Math., Vol 87(2002), 5-36.


\bibitem{Bo} R. Bowen. {\em Hausdorff dimension of quasicircles.} Publ. Math. IHES Vol50(1979), 11-25.

 \bibitem{C}L. Carleson. {\em On $H^{\infty}$ in multiply connected domains},Conference in harmonic analysis in honor of Antoni Zygmund, Wadsworth, 1983.

\bibitem{Cu}G. Cui, {\em Integrably asymptotic affine homeomorphisms of the circle and {\teich} spaces}, Sci. China Ser A, 43, (2000), 267-279.
%\bibitem{Du} P. Duren, {\em theorem of $H^{p}$ spaces}, New York: Academic Press, 1970.
\bibitem{Da}F. Dal'Bo. {\em Geodesic and horocyclic Trajectories}£¬ Springer, 2011.
\bibitem {F} J.L. Fernandez, {\em Domains with strong barrier,} Revista. Math.. Iber. 5, (1989), 47-65.
\bibitem {FH} J.L. Fernandez, D.H. Hamilton, {\em Length of curves under conformal mappings,} Comm. Math. Helv. 62, (1987), 122-134.
\bibitem{Ga}J. B. Garnett: {\em  Bounded Analytic Functions,} New York: Academic Press,1981.

\bibitem{H1} S. Huo:  {\em On Carleson measures induced by Beltrami coefficients being compatible with Fuchsian groups,} Ann. Acad. Sci. Fenn. Math., Vol 46(2021).
\bibitem{HW} S. Huo and S. Wu:  {\em The failure of analyticity of Hausdorff dimensions of quasi-circles of Fuchsian groups of the second kind,} Proc. Amer. Math. Soc. Vol 143(2015), 1101-1108.

\bibitem {HZ} S. Huo and M. Zinsmerster: {\em On Ruelle's property}. arXiv:1906.01291v1.
\bibitem {KMS} Y. Komori, V. Markovic and C. Series: {\em Kleinian groups and hyperbolic 3-Manifolds}, London Math. Soc. Lecture note series,299, 2003.


\bibitem{Po} C. Pommerenke, {\em Schlichte Functionen und BMOA }, Comment. Helv.,52, (1977), 591--602.

\bibitem{RR}L. A. Rubel and J.V. Ryff, {\em The bounded weak-star topology and the bounded analytic functions. } J. Funct. Anal.,Vol 5(1970), 167-183.

\bibitem[11]{Ru} D. Ruelle. {\em Repellers for real analytic maps.} Ergod. Th. and Dynam.Syst., 2(1982)99-107.

\bibitem{SW}Y. Shen and H. Wei, {\em Universal Teichmuller space and BMO}, Adv. Math, Vol 234(2013), 129-148.


\bibitem{Se} S. Semmes, {\em Quasiconformal mappings and chord-arc curves}. Trans. Am. Math. Soc. 306(1988),233-263.

%
\bibitem{Su1} D. Sullivan, {\em Growth of positive harmonic functions and kleinian group limit sets of zero plane measure and Hausdorff dimension two}, Geometry Symposium Vol 504. (1980), 127-144.
%
\bibitem{Su2} D. Sullivan, {\em Discrete conformal groups and measurable dynamics}, Bull. Amer. Math. Soc. Vol 6. (1982), 57-73.

%%
\bibitem {WZ}Wei H. and M. Zinsmeister: {\em Carleson  measures and chord-arc curves}, Ann. Acad. Sci. Fenn. Math., Vol 43 (2018),466-483.

%
\end{thebibliography}
\end{document}